\documentclass{amsart}
\usepackage[T1]{fontenc}
\usepackage{ae}
\usepackage{aecompl}%
\usepackage{ucs}
\usepackage[utf8x]{inputenc}
\usepackage{color}
\usepackage[enableskew,vcentermath]{youngtab}
\usepackage{amsfonts}
\usepackage{hhline}
\usepackage{amsthm}
\usepackage{amsmath}
\usepackage{amssymb}
\usepackage[all]{pstricks}
\usepackage{hyperref}

\begin{document}

\numberwithin{equation}{section}

\newtheorem{Le}{Lemma}[section]
\newtheorem{Ko}[Le]{Lemma}
\newtheorem{Sa}[Le]{Theorem}
\newtheorem{pro}[Le]{Proposition}

\newtheorem{theorem}[Le]{Theorem}
\newtheorem{lemma}[Le]{Lemma}
\newtheorem{Con}[Le]{Conjecture}
\newtheoremstyle{Bemerkung}
  {}{}{}{}{\bfseries}{.}{0.5em}{{\thmname{#1}\thmnumber{ #2}\thmnote{ (#3)}}}

\theoremstyle{Bemerkung}
\newtheorem{definition}[Le]{Definition}
\newtheorem{Def}[Le]{Definition}
\newtheorem{example}[Le]{Example}

\newtheorem{remark}[Le]{Remark}
\newtheorem{Bem}[Le]{Remark}
\newtheorem{Bsp}[Le]{Example}

\renewcommand{\l}{\lambda}
\renewcommand{\L}{\lambda}
\newcommand{\bl}{\bar\lambda}
\newcommand{\bn}{\bar\nu}
\newcommand{\A}{\mathcal{A}}
\newcommand{\B}{\mathcal{B}}
\newcommand{\C}{\mathcal{C}}
\newcommand{\D}{\mathcal{D}}
\renewcommand{\S}{\mathcal{S}}
\renewcommand{\L}{\mathcal{L}}
\newcommand{\E}{\mathcal{E}}
\newcommand{\N}{\mathcal{N}}
\renewcommand{\a}{\alpha}
\renewcommand{\b}{\beta}
\renewcommand{\c}{\gamma}
\renewcommand{\d}{\delta}
\newcommand{\e}{\epsilon}
\newcommand{\h}{\hfil}
\newcommand{\X}{X}
\newcommand{\abs}[1]{\left| #1 \right|}
\newcommand{\lm}{\l/\mu}
\renewcommand{\ln}{\l/\nu}
\newcommand{\ab}{\a/\b}
\newcommand{\m}{\mu}
\newcommand{\n}{\nu}

\newcommand{\LR}{\mathcal{LR}}
\newcommand{\LRk}{\mathcal{LR}_k}

\newcommand{\Rk}{\mathbb{R}^{3k}}

\title[Generalised Stretched LR Coefficients]{Generalised Stretched Littlewood-Richardson Coefficients}
\author[C. Gutschwager]{Christian Gutschwager}
\address{Institut für Algebra, Zahlentheorie und Diskrete Mathematik, Leibniz Universität Hannover,  Welfengarten 1, D-30167 Hannover}
\email{gutschwager (at) math (dot) uni-hannover (dot) de}

\subjclass[2000]{05E05,05E10,14M15,20C30}
\keywords{Littlewood Richardson, LR coefficient, LR tableaux}

\begin{abstract}
The Littlewood-Richardson (LR) coefficient counts, among many other things, the LR tableaux of a given shape and a given content. We prove that the number of LR tableaux weakly increases if one adds to its shape and content the shape and the content of another LR tableau. We also investigate the behaviour of the number of LR tableaux, if one repeatedly adds to the shape another shape with either fixed or arbitrary content. This is a generalisation of the stretched LR coefficients, where one repeatedly adds the same shape and content to itself.
\end{abstract}

\maketitle

\section{Introduction}
The Littlewood-Richardson (LR) coefficients $c(\l;\m,\n)$ appear in many branches of mathematics. For example, they appear  in the representation theory of the symmetric groups, in the theory of symmetric functions, in the Schubert calculus and in problems regarding the existence of matrices with certain eigenvalues or invariant factors (see \cite{Fulton}).

Some recent research has been concerned with the behaviour of the stretched LR coefficients. More precisely, fix partitions $\l,\m,\n$ and investigate the function $f(n)=c(n\l;n\m,n\n)$ as a function of $n$, where $n\l$ is the partition obtained from $\l$ by multiplying every part by $n$. King et al. \cite{S3} conjectured that the stretched LR coefficient is a polynomial in $n$. Derksen and Weyman \cite{DW} as well as Rassart \cite{Rassart} proved King's conjecture to hold true, using semi-invariants of quivers and partition functions, respectively. In fact, Rassart \cite{Rassart} proved even more: fix a positive integer $k$ and let the partitions $\l,\m,\n$ have length at most $k$, then the triples $(\l,\m,\n)$ of partitions with positive LR coefficient $c(\l;\m,\n)$ form a cone in $\Rk$. This cone decays into a finite number of cones in which the LR coefficient is given by a polynomial in $(\l_1,\l_2,\ldots,\m_1,\ldots,\n_k)$. Rassart remarks in his paper that Knutson also has an unpublished proof for this property using symplectic geometric techniques.

 In  \cite{KnutsonTao}  (see also \cite{S1}) Knutson and Tao proved the saturation conjecture which is, that $f(n)=c(n\l;n\m,n\n)\neq 0$ for some $n\geq 1$ implies $c(\l;\m,\n)\neq 0$.
In  \cite{KTW} Knutson, Tao and Woodward  proved that $f(n)$ is constant if and only if $c(\l;\m,\n)=1$. Furthermore, if the polynomial $f(n)\neq0$ has an integer root $-t\in\mathbb{Z}$ then $t>0$ and  $f(n)$ also contains the factors $(n+i)$ for $1\leq i\leq t$. Furthermore, there is a $t$ such that  $f(n)=g(n)\prod_{i=1}^t(n+i)$ with $g(n)$ a polynomial with no integer roots.  Let $\l+\l'=(\l_1+\l'_1,\l_2+\l'_2,\ldots)$. We will show in Lemma~\ref{Le:add} that $c(\l';\m',\n')\neq 0$ implies $ c(\l+\l';\m+\m',\n+\n')\geq c(\l;\m,\n)$.

Our main object of study   is an affine generalisation of the stretched LR coefficient, namely $P(n)=P^{\l,\m,\n}_{\l',\m',\n'}(n)=c(n\l+\l';n\m+\m',n\n+\n')$. Using Lemma~\ref{Le:add} we will in Remark~\ref{Bem33} make a first observation about $P(n)$, namely that $P(n)$ is weakly increasing. 
To obtain more results about $P(n)$, in Section~\ref{seq:Q} we will  investigate the function $Q(n)=Q^{\l,\m}_{\l',\m'}(n)=\sum_\n c(n\l+\l';n\m+\m',\n)$ in more detail. The function $Q(n)$ counts the LR tableaux of shape $n\l+\l'/n\m+\m'$ and arbitrary content, which is therefore the number of irreducible characters (counted with multiplicity) in the skew character $[n\l+\l'/n\m+\m']$. Our main result will be that $Q(n)$ is bounded above if   and only if  $\lm$ is a partition or a rotated partition (see Theorem~\ref{Sa:main} for the if part and Lemma~\ref{Q:prop} for the only if part). Furthermore, if $\lm$ is a partition or a rotated partition, then $Q(n)$ is strictly increasing until it reaches its upper bound. In Theorem~\ref{Sa:main} we also give the value $n$ for which $Q(n)$ at first obtains the upper bound.

In Section~\ref{sec:P}, we will  investigate the generalised stretched LR coefficient $P(n)=c(n\l+\l';n\m+\m',n\n+\n')$ as a function of $n$ in more detail. We will see in Lemma~\ref{Le:P} that $P(n)$ has an upper bound in some cases, for example if $\lm$ is a partition or a rotated partition. Furthermore, by Lemma~\ref{Le:Rassartg}, for large $n$ the function $P(n)$ is given by a polynomial, which has, by Lemma~\ref{Le:P2}, in some cases the same degree as the polynomial $c(n\l;n\m,n\n)$.

\section{Notation and Littlewood-Richardson Symmetries}
We mostly follow the standard notation in \cite{Sag} or \cite{Stanley}. A partition $\l=(\l_1,\l_2,\ldots,\l_l)$ is a weakly decreasing sequence of non-negative integers where only finitely many of the $\l_i$ are positive. We regard two partitions as the same if they differ only by the number of trailing zeros and call the positive $\l_i$ the parts of $\l$. The length is the number of positive parts and we write $l(\l)=l$ for the length and $\abs{\l}=\sum_i \l_i$ for the sum of the parts. With a partition $\l$ we associate a diagram, which we also denote by $\l$, containing $\l_i$ left-justified boxes in the $i$-th row and we use matrix style coordinates to refer to the boxes.

The {conjugate} $\l^c$ of $\l$ is the diagram which has $\l_i$ boxes in the $i$-th column.

The sum $\m+\n=\l$ of two partitions $\m,\n$ is defined by $\l_i=\m_i+\n_i$. The partition $\m\cup\n$ contains the parts of both $\m$ and $\n$. These operations are conjugate to another
\[(\m+\n)^c=\m^c\cup\n^c.\]
For example, we have
\begin{align*}
 \yng(4,4,2,1)+\young(XXX,XX,X,X)&=\young(\h\X\h\h\h\X\X,\h\X\h\h\h\X,\h\X\h,\h\X), & \yng(4,3,2,2) \cup \young(XXXX,XX,X)&=\young(\h\h\h\h,XXXX,\h\h\h,\h\h,\h\h,XX,X).
\end{align*}

For $\m \subseteq \l$ we define the skew diagram $\lm$ as the difference of the diagrams $\l$ and $\m$, defined as the difference of the set of the boxes. Rotation of $\lm$ by $180^\circ$ yields a skew diagram $(\lm)^\circ$ which is well defined up to translation. A skew tableau $T$ is a skew diagram in which positive integers are written into the boxes. A semistandard tableau of shape $\lm$ is a filling of $\lm$ with positive integers such that the entries weakly increase along the rows and strictly increase down the columns. The content of a semistandard tableau $T$ is $\n=(\n_1,\ldots)$ if the number of occurrences of the entry $i$ in $T$ is $\n_i$. The reverse row word of a tableau $T$ is the sequence obtained by reading the entries of $T$ from right to left and top to bottom starting with the first row. Such a sequence is said to be a lattice word if for all $i,n \geq1$ the number of occurrences of $i$ among the first $n$ terms is at least the number of occurrences of $i+1$ among these terms. The Littlewood-Richardson (LR) coefficient $c(\l;\m,\n)$ equals the number of semistandard tableaux of shape $\lm$ with content $\n$ such that the reverse row word is a lattice word. We will call those tableaux LR tableaux. The LR coefficients play an important role in different contexts (see \cite{Fulton,Sag,Stanley} for further details).

The irreducible characters $[\l]$ of the symmetric group $S_n$ are naturally labeled by partitions $\l\vdash n$. The skew character $[\lm]$ corresponding to a skew diagram $\lm$ is defined by the LR coefficients
\[ [\lm]=\sum_\n c(\l;\m,\n) [\n]. \]

Let $\A$ and $\B$ be non-empty subdiagrams of a skew diagram $\D$ such that the union of $\A$ and $\B$ is $\D$.
Then we say that the skew diagram $\D$ is \textit{disconnected} or \textit{decays} into the skew diagrams $\A$ and $\B$ if no box of $\A$ (viewed as boxes in $\D$) is in the same row or column as a box of $\B$. We write $\D=\A\otimes\B$ if $\D$ decays into $\A$ and $\B$. A skew diagram is \textit{connected} if it does not decay. If $\D=\A\otimes\B=\C$ then by translation symmetry $[\D]=[\C]$.

For example, the skew diagram $\D=\young(:::::\h\h\h,:::\h\h\h\h\h,:::\h,::\h,:\h\h,\h,\h)$ is disconnected and decays into the skew diagrams $(5,5,1)/(2), (2,2)/(1)$ and $(1^2)$ which are connected. So we have $\D=(5,5,1)/(2)\otimes(2,2)/(1)\otimes(1^2)$.

\textit{Translation symmetry} gives $[\lm]=[\ab]$ if the skew diagrams of $\lm$ and $\ab$ are the same up to translation. Translation includes the case that $\lm$ decays and  connected subdiagrams are translated independent of each other.   Furthermore, \textit{rotation symmetry} gives $[(\lm)^\circ]=[\lm]$.  \textit{Conjugation symmetry} $c(\l^c;\m^c,\n^c)=c(\l;\m,\n)$ is also well known, as is $c(\l;\m,\n)=c(\l;\n,\m)$.

A \textit{basic skew diagram} $\lm$ is a skew diagram which satisfies $\m_i<\l_i$ and $\m_i\leq \l_{i+1}$ for each $1\leq i \leq l(\l)$. This means that we do not have empty rows or columns in $\lm$. Empty rows or columns of a skew diagram do not influence the filling and so deleting empty rows or columns does not change the skew character or LR fillings.

A \textit{proper skew diagram} $\lm$ is a skew diagram which is neither a partition nor a rotated partition.

In~\cite{Gut} we used the following theorem to classify  multiplicity free skew characters.
\begin{Sa}[Theorem 3.1, \cite{Gut}]
\label{Sa:addold}
Let $\l, \m, \n $ be partitions and $a,b\geq 0$ be integers. Then
\[c(\l+(1^{a+b});\m+(1^a),\n+(1^{b}))\geq c(\l;\m,\n) \]
and by conjugation
\[ c(\l\cup(a+b);\m\cup(a),\n\cup(b))\geq c(\l;\m,\n)  .\]

\end{Sa}

\section{Key Lemma}

We can generalise Theorem~\ref{Sa:addold} to the following.
\begin{Le}
\label{Le:add}
Let $\l,\m,\n,\l',\m',\n' $ be partitions with $c(\l;\m,\n),c(\l';\m',\n')\neq0$.
Then
\[  c(\l+\l';\m+\m',\n+\n')\geq c(\l;\m,\n)\]
and by conjugation
\[ c(\l\cup\l';\m\cup\m',\n\cup\n')\geq c(\l;\m,\n).\]
\end{Le}
\begin{proof}
Let $\A$ be a fixed LR tableau of shape $\l'/\m'$ with content $\n'$. Let $A_j$ be the multiset of the entries in the $j$th row of $\A$.

For any LR tableau $\C^i$ of shape $\l/\m$ and content $\n$ we let $\Gamma^i_j$ be the multiset of the entries in the $j$th row of $\C^i$.

We can now define for every $\C^i$ a tableau $\D^i$ of shape $(\l+\l')/(\m+\m')$ with content $\n+\n'$ by placing  the entries of $A_j\cup \Gamma^i_j$ into row $j$ in weakly increasing order. To see that the entries are strictly increasing down the columns let $\C^i_j$ denote the multiset of the entries of the $j$th row of $\C^i$ where we assume that the empty boxes belonging to $\m$ contribute a $0$ each. So there are $\m_j$ additional entries $0$ in $\C^i_j$ compared to $\Gamma^i_j$. Define $\D^i_j$ and $\A_j$ accordingly. Clearly $\D^i_j=\A_j\cup\C^i_j$. Now the entries in $\C^i$ are increasing down the columns if and only if the number of entries smaller than or equal to $k$ in $\C^i_j$ is at most the number of entries smaller than $k$ in $\C^i_{j-1}$ for each $k,j>1$. Since $\A$ and $\C^i$ are semistandard, $\C^i_j$ and $\A_j$ satisfy this condition and so does $\D^i_j$, hence $\D^i$ is semistandard.

It is also clear that the tableau word is a lattice word because it can be divided into two subsequences (corresponding to the entries in $\D^i$ having their origin in either $\A$ or $\C^i$) which are both lattice words. So the $\D^i$ are in fact LR tableaux.

Suppose we have $\D^i=\D^l$. Then we know by construction that the multiset of the entries in the $j$th row of $\D^i$ is $A_j\cup \Gamma^i_j$ while the multiset of the entries in the $j$th row of $\D^l$ is $A_j\cup \Gamma^l_j$. This gives us $\Gamma^i_j=\Gamma^l_j$ for all $j$ and, since an LR tableau of a given shape is uniquely determined by the content of its rows, it follows that $\C^i=\C^l$. So different LR tableaux of shape $\lm$ with content $\n$ give different LR tableaux of shape $(\l+\l')/(\m+\m')$ with content $\n+\n'$, thus
\[c(\l;\m,\n) \leq c(\l+\l';\m+\m',\n+\n').\qedhere \]
\end{proof}
\begin{Bem}\label{Bem:king}
 In the hive model (which we do not use in this paper) the proof is also easy. Choose one LR hive corresponding to the triple $(\l',\m',\n')$ and add this hive to all the LR hives corresponding to $(\l,\m,\n)$. It is easy to see that all the new hives are different LR hives corresponding to $(\l+\l',\m+\m',\n+\n')$.
\end{Bem}

\begin{Bem}\label{Bem33}
It is known that $f(n)=c(n\l;n\m,n\n)$ is a polynomial which is constant if and only if $c(\l;\m,\n)=1$ (see \cite{KnutsonTao},\cite{KTW}). Suppose $\l,\m,\n$ are chosen in such a way that $f(n)$ is not constant. Then it follows from Lemma~\ref{Le:add} that
\begin{align*}
P(n)&=c(n\l+\l';n\m+\m',n\n+\n')
\end{align*}
increases without bound if $c(\l';\m',\n')\neq 0$.
\end{Bem}

\begin{Bem}
It is known (see \cite{Zelevinsky}) that the triples of partitions with non-zero LR coefficient form an additive semigroup.
\end{Bem}

\begin{Le}\label{Le:cons}
 Let $\lm$ and $\l'/\m'$ be skew diagrams. Then
\[\sum_{\n'} c(\l+\l';\m+\m',\n') \geq \sum_\n c(\l;\m,\n). \]
\end{Le}
\begin{proof}
 Since $\l'/\m'$ is a skew diagram there exists a partition $\a$ which satisfies $c(\l';\m',\a)\neq0$. Note that if $\l'/\m'$ is empty we can choose $\a=\emptyset$ and have $c(\l';\m',\a)=1$. Now different $\n$ give different $\n+\a$ and by Lemma~\ref{Le:add} we have 
 \[c(\l+\l';\m+\m',\n+\a) \geq  c(\l;\m,\n). \]

So $\sum_\n c(\l+\l';\m+\m',\n+\a) \geq \sum_\n c(\l;\m,\n)$.  Extending the sum on the left hand side from $\n+\a$ to arbitrary $\n'$ gives

\[\sum_{\n'} c(\l+\l';\m+\m',\n') \geq \sum_\n c(\l+\l';\m+\m',\n+\a) \geq\sum_\n c(\l;\m,\n).\qedhere\]
\end{proof}

\section{The Behaviour of $Q^{\l,\m}_{\l',\m'}(n)$}\label{seq:Q}

For $\m \subseteq \l,\m' \subseteq \l'$ we define $Q^{\l,\m}_{\l',\m'}(n)=\sum_\n c(n\l+\l';n\m+\m',\n)$ and write simply $Q(n)$ if $\l,\m,\l',\m'$ are known from the context.

\begin{Le} \label{Q:prop}
 Let $\lm$ be a proper skew diagram. Then $Q^{\l,\m}_{\l',\m'}(n)$ increases without bound as $n$ increases.
Furthermore,
\begin{align*}
\abs{\{ \n |  c(n\l+\l';n\m+\m',\n)\neq 0 \}  } &\rightarrow \infty & \textnormal{as } n &\rightarrow \infty.
\end{align*}

\end{Le}
\begin{proof}
Since $\lm$ is a proper skew diagram, it is obtained from the skew diagram $(2,1)/(1)$ by inserting rows and columns and so by Lemma~\ref{Le:cons}  we have $\sum_\n c(\l;\m,\n) \geq \sum_\n c( (2,1);(1),\n) $. Furthermore, we have
\[ \sum_\n c(n\l+\l';n\m+\m',\n) \geq  \sum_\n c(n\l;n\m,\n) \geq \sum_\n c( n(2,1);n(1),\n) .\]
It is easy to see that $\sum_\n c( n(2,1);n(1),\n)=n+1$, because an LR tableau of shape $(2n,n)/(n)$ contains $n$ entries $1$ in the first row and $i$ ($0\leq i \leq n$) entries $1$ as well as $n-i$ entries $2$ in the second row. Moreover, for each such $i$ there is exactly one LR tableau of shape $(2n,n)/(n)$. Thus,  $Q^{\l,\m}_{\l',\m'}(n)$ increases without bound.

Furthermore, since the number of irreducible characters in $[n(2,1)/n(1)]$ is $n+1$, there are also at least $n+1$ irreducible characters in $[n\l+\l'/n\m+\m']$ (by the same argument as in Lemma~\ref{Le:cons}), hence
\[\abs{\{ \n |  c(n\l+\l';n\m+\m',\n)\neq 0 \}  }= \sum_{\genfrac{}{}{0pt}{}{\n} {c(n\l+\l';n\m+\m',\n)\neq 0}} 1 \geq n+1 .\qedhere\]
\end{proof}

\begin{Sa}\label{Sa:main}
Let $\lm$ be a partition or a rotated partition.
Then there exists an $m$ with  $Q^{\l,\m}_{\l',\m'}(n)=Q^{\l,\m}_{\l',\m'}(m)$ for $n \geq m$.
Furthermore, suppose $\l=(\a_1^{a_1},\a_2,\a_3,\ldots\a_k),\a_k\neq 0,\: \m=(\a_1^{a_1-1})$ and $\l'/\m'$ basic. Then the smallest $m$ we can choose for the above equation is given by
\[ m=\biggl\lceil \max_{ \genfrac{}{}{0pt}{}{1\leq j \leq k} {\a_j>\a_{j+1}}} \Bigl( \frac{\l'_1- \l'_{a_j}+ \l'_{a_j+1}+\m'_{a_1}-\m'_{a_1-1}}{\a_j-\a_{j+1}} \Bigr)   \biggr\rceil\]
(where $\lceil x \rceil$ denotes as usual the smallest integer larger than  or equal to $x$) with $a_j= a_1-1+j,\a_{k+1}=0$ (for $a_1=1$ set $\m'_{0}=\l'_1$). Furthermore, 
\[Q^{\l,\m}_{\l',\m'}(m)>Q^{\l,\m}_{\l',\m'}(m-1)>\ldots>Q^{\l,\m}_{\l',\m'}(0).\]

These inequalities are also satisfied in the general case if we choose the smallest $m$ satisfying $Q^{\l,\m}_{\l',\m'}(n)=Q^{\l,\m}_{\l',\m'}(m)$ for $n \geq m$.
\end{Sa}
\begin{proof}
We look at the skew diagram $\A(n)=(n\l+\l')/(n\m+\m')$.

By rotation symmetry we may assume that $\lm$ is a partition instead of a rotated partition.

Let $a_1>a_2>\ldots>a_k$ be the indices of the non-empty rows of $\lm$. If we have $\l_i=\m_i>\l_{i+1}$ for some $i\neq a_1,...,a_k$ and choose $n$ big enough then $\A(n)$ decays into a skew diagram $\A^{up}$ containing the top $i$ rows and a skew diagram $\A_{lo}$ containing the rows below row $i$. If we increase $n$ even more then the skew diagrams $\A^{up}$ and $\A_{lo}$ are translated relative to one another which is irrelevant for the skew character $[\A(n)]$. So if there are some $i\neq a_1,...,a_k$ with $\l_i=\m_i>\l_{i+1}$ we may choose $n$ large enough so that for each such $i$, $\A(n)$ decays into an upper skew diagram and a lower skew diagram. Instead of looking at this situation we may then investigate the case that $\l'/\m'=\A(n)$ for an $n$ large enough and have no $i\neq a_1,...,a_k$ with $\l_i=\m_i>\l_{i+1}$. So we may assume that $\m_i=\l_i=\l_{a_1}$ for $i<a_1$ and $\m_i=\l_i=\m_{a_k}$ for $a_k<i\leq l(\m)$ (and since $\l/\m$ is a partition we also have $\m_{a_1}=\m_{a_k}$). If $\m_{a_1}>0$ there is for the same reason as above an $n$ such that $\A(n)$ decays into skew diagrams containing the top $l(\m)$ rows and the rows below row $l(\m)$ and increasing $n$ translates these skew diagrams relative to another so we may assume that $\m_{a_1}=0$.

As an example for the above, assume
\begin{figure}[h]
\psset{xunit=0.4cm,yunit=0.4cm,runit=0.4cm}
\begin{pspicture*}(-0.2,2,8)(15.2,6.2)
\put(0,4.5){$\lm=(5,4,2)/(5,2,2)=$} \psline (10,3)(12,3)(12,4)(14,4)(14,5)(15,5)(15,6)
\psline(12,4)(12,5)(14,5)
\psline(13,4)(13,5)  \psdot(11,3)\psdot(12,3)\psdot(15,5)
\end{pspicture*}
\end{figure}
\newline and $\A(0)=\l'/\m'=(5,4,3,3)/(2,1)=\young(::\h\h\h,:\h\h\h,\h\h\h,\h\h\h)$.
We then have
\begin{align*}
 \A(1)&=\young(:::::::\h\h\h,:::\h\h\h\h\h,::\h\h\h,\h\h\h)\\
 \A(2)&=\young(::::::::::::\h\h\h,:::::\h\h\h\h\h\h\h,::::\h\h\h,\h\h\h)\\
 \A(3)&=\young(:::::::::::::::::\h\h\h,:::::::\h\h\h\h\h\h\h\h\h,::::::\h\h\h,\h\h\h).
\end{align*}
So for $n\geq2$ the skew diagram $\A(n)$ decays into three connected skew diagrams. and the only effect of the empty columns of $\lm$ for $n\geq2$ is that those three skew diagrams are translated relative to another. But since translation is irrelevant for LR fillings we can instead investigate the situation $\lm=(2,2)/(2)$ and $\l'/\m'=\A(2)=(14,11,6,3)/(11,4,3)=\young(:::::::::::\h\h\h,::::\h\h\h\h\h\h\h,:::\h\h\h,\h\h\h)$ where we additionally removed the empty column to make $\l'/\m'$ basic.

So we have, without loss of generality, $\l=(\a_1^{a_1},\a_2,\a_3,\ldots\a_k),\a_k\neq0$ (not necessarily $\a_i\neq\a_{i+1}$) and $\m=(\a_1^{a_1-1})$. To prove $Q(n)=Q(m)$ for $n \geq m$, we have to construct an $m$ such that removing in an LR tableau of shape $\A(n)$ from the row $a_i$ with $1\leq i \leq k$ the entry $i$ $(n-m)\a_i$ times and translating the top $a_1-1$ rows $(n-m)\a_1$ boxes to the left yields an LR tableau of shape $\A(m)$.

By our choice of $\l$ and $\m$, the number $N$ of non-empty columns among the top $a_1-1$ rows of $\A(n)$ is independent of $n$. We have $N\leq \l'_1-\m'_{a_1-1}$ and may by translation symmetry assume equality (set $\m'_0=\l'_1$ for $a_1=1$). So the number of entries $1$ among the top $a_1-1$ rows of an LR filling of $\A(n)$ is at most $N$. So for $1\leq i \leq k$ there are at most $N$ entries larger than $i$ in row $a_i$ of an LR filling of $\A(n)$. Furthermore, the number of entries smaller than $i$ in row $a_i$ is at most $\m'_{a_1}-\m'_{a_i}$, which is also independent of $n$. On the other hand, there are $\l'_{a_i}-\m'_{a_i}+n\a_i$ boxes in row $a_i$ of $\A(n)$. So the number of entries $i$ in row $a_i$ of an LR filling of $\A(n)$ is at least \[\l'_{a_i}-\m'_{a_i}+n\a_i - N - (\m'_{a_1}-\m'_{a_i}) = \l'_{a_i}-\m'_{a_1}- N +n\a_i.\] Obviously, if $\l'_{a_k}-\m'_{a_1}- N +n\a_k \geq 0$ then also $\l'_{a_i}-\m'_{a_1}- N +n\a_i \geq 0$ for every $1\leq i \leq k$. So for
\begin{equation} \label{eqQ'}
n>n'\geq\frac{\m'_{a_1}+N-\l'_{a_k}}{\a_k}
\end{equation}
there are at least $(n-n')\a_i$ entries $i$ in row $a_i$ of every LR tableau of shape $\A(n)$.

We have to investigate the $j$ ($1\leq j \leq k$) with $\a_j>\a_{j+1}$ (for example $j=k$). Removing $\a_i$ times the entry $i$ from row $a_i$ in an LR tableau removes more entries $j$ than $j+1$ so the new tableau can violate the lattice word condition even if there are enough entries $i$ to remove. As calculated above the number of entries $j$ in row $a_j$ of an LR tableau of shape $\A(n)$ is at least $\l'_{a_j}-\m'_{a_1}- N +n\a_j$. Furthermore, the number of entries $j+1$ below row $a_j$ in an LR tableau of shape $\A(n)$ is at most $\l'_{a_j+1}+n\a_{j+1}$ since this is the number of columns below row $a_j$. So for \[\l'_{a_j}-\m'_{a_1}- N +n\a_j \geq \l'_{a_j+1}+n\a_{j+1}\] the number of entries $j$ in row $a_j$ is at least as large as the number of entries $j+1$ below row $a_j$ in every LR tableau of shape $\A(n)$. We can solve the above inequality and get
\[n \geq \frac{\l'_{a_j+1}- \l'_{a_j}+\m'_{a_1}+ N}{\a_j-\a_{j+1}}.\]

Since we have $\a_k>0=\a_{k+1}$ setting $j=k$ gives \[\frac{ \l'_{a_k+1}-\l'_{a_k}+\m'_{a_1}+N}{\a_k}\geq \frac{-\l'_{a_k}+\m'_{a_1}+N}{\a_k} \] which is the right hand side of inequality~\eqref{eqQ'}.

Let us set
\[ m=\biggl\lceil \max_{ \genfrac{}{}{0pt}{}{1\leq j \leq k} {\a_j>\a_{j+1}}} \Bigl( \frac{\l'_1- \l'_{a_j}+ \l'_{a_j+1}+\m'_{a_1}-\m'_{a_1-1}}{\a_j-\a_{j+1}} \Bigr)   \biggr\rceil.\]

Then we know from the arguments above, that for $n\geq m$ every LR tableau $\C_n$ of shape $\A(n)$ contains at least $(n-m)\a_i$ entries $i$ in row $a_i$ ($1\leq i\leq k$). Furthermore, removing $(n-m)\a_i$ entries $i$ from every row $a_i$ ($1\leq i\leq k$) and translating the top $a_1-1$ rows $(n-m)\a_1$ boxes to the left yields a tableau $\C_m$ which contains (for those $j$ with $\a_j>\a_{j+1}$) at least as many entries $j$ in row $a_j$  as there are entries $j+1$ below row $a_j$. So the tableau $\C_m$ satisfies the lattice word condition. Furthermore, the entries in the rows weakly increase  from left to right. We have to check that the entries in the columns are strictly increasing from top to bottom which is non trivial because we remove more entries $j$ from row $a_j$ than entries $j+1$ from row $a_j+1$ if $\a_j>\a_{j+1}$. The condition on $m$ ensures that in $\C_m$ there is an entry smaller than $j+1$ above every entry in row $a_j+1$ so there is no problem for the entries greater than or equal to  $j+1$ in row $a_j+1$. Furthermore, the entries in $\C_m$ in row $a_j+1$ which are smaller than $j+1$ have an entry smaller than itself in the box directly above itself because $\C_n$ is semistandard. So $\C_m$ is in fact an LR tableau. So every LR tableau of shape $\A(n)$ is obtained from an LR tableau of shape $\A(m)$ by adding $(n-m)\a_i$ entries to row $a_i$ ($1\leq i \leq k$) and translating the top $a_1-1$ rows $(n-m)\a_1$ boxes to the right. So for $n\geq m$ we have $Q(n)=Q(m)$.

We now have to prove that  $Q^{\l,\m}_{\l',\m'}(m)>Q^{\l,\m}_{\l',\m'}(m-1)>\ldots>Q^{\l,\m}_{\l',\m'}(0)$ if $\l'/\m'$ is basic.

For $n<\frac{\l'_1-\l'_{a_k}+\m'_{a_1}-\m'_{a_1-1}}{\a_k}$ we can construct an LR tableau of shape $\A(n)$ containing fewer than $\a_k$ entries $k$ in row $a_k$ (see the example below). The existence of such an LR tableau follows directly from the  arguments above and gives $Q(n)>Q(n-1)$.

Now suppose $\frac{\l'_1-\l'_{a_k}+\m'_{a_1}-\m'_{a_1-1}}{\a_k}\leq n<\frac{\l'_1- \l'_{a_j}+ \l'_{a_j+1}+\m'_{a_1}-\m'_{a_1-1}}{\a_j-\a_{j+1}}$ for some $1\leq j \leq k$ with $\a_j>\a_{j+1}$. We can construct an LR tableau $C_n$ of shape $\A(n)$ satisfying the following conditions (also see the examples below).
\begin{itemize}
 \item There are $\l'_1-\m'_{a_1-1}$ entries $1$ in the top $a_1-1$ rows of $C_n$ (this is possible because $\l'/\m'$ is basic).
 \item There are $\l'_1-\m'_{a_1-1}$ entries $2$ in the top $a_2-1$ rows of $C_n$ (the lower bound on $n$ ensures that there are enough boxes in row $a_2-1$).
\item For $1\leq i \leq j$, there are $\l'_1-\m'_{a_i-1}$ entries $i$ in the top $a_i-1$ rows of $C_n$ (the lower bound on $n$ ensures that there are enough boxes in row $a_i-1$).
 \item There are $\l'_1-\m'_{a_1-1}$ entries $j$ in the top $a_j-1$ rows of $C_n$ (the lower bound on $n$ ensures that there are enough boxes in row $a_j-1$).
 \item There are $\l'_1-\m'_{a_1-1}$ entries $j+1$ in the top  $a_j$ rows (the lower bound on $n$ ensures that there are enough boxes in row $a_j$).
 \item There are at least $x\geq \a_j$ entries $1$ in row $a_1$. For $2\leq i <j$ there are at least $x$ entries $i$ in row $a_i$ and there are exactly $x$ entries $j$ in row $a_j$ and $x$ entries $j+1$ below row $a_j$ (the upper bound on $n$ ensures that there are at least $x$ columns below row $a_j$ into which we can write the entry $j+1$).
 \item There is no entry $j$ below row $a_j$.
 \item Fill the other boxes, for example, in increasing order for each column.
\end{itemize}
 By construction $C_n$ is an LR tableau. Removing  $\a_i$ entries $i$ from every row $a_i$ and translating   the top $a_1-1$ rows by $\a_1$ boxes to the left, yields a tableau $C_{n-1}$ which contains more entries $j+1$ than entries $j$ and so is not an LR tableau. This gives $Q(n)>Q(n-1)$.

Take for example
\begin{figure}[h]
\psset{xunit=0.4cm,yunit=0.4cm,runit=0.4cm}
\begin{pspicture*}(-0.2,2,8)(14.2,7.2)
\put(0,5){$\lm=(3,3,3,1)/(3,3)=$} 
\psline(12,5)(12,3)(11,3)(11,5)(14,5)(14,4)(11,4)
\psline(13,5)(13,4)
\psline(14,5)(14,7)\psdot(14,6)
\end{pspicture*}
\end{figure}

\[\l'/\m'=(12,11,10,9,5,3,3,1)/(8,6,6,3,1,1,1)=\young(::::::::\h\h\h\h,::::::\h\h\h\h\h,::::::\h\h\h\h,:::\h\h\h\h\h\h,:\h\h\h\h,:\h\h,:\h\h,\h).\]
We now want to construct the aforementioned LR tableaux. We have $a_1=3,a_2=4,\a_1=3,\a_2=1,\a_3=0$ and $k=2$ and therefore
\begin{itemize}
 \item $\frac{\l'_1-\l'_{a_k}+\m'_{a_1}-\m'_{a_1-1}}{\a_k}=3$,
\item $\frac{\l'_1- \l'_{a_j}+ \l'_{a_j+1}+\m'_{a_1}-\m'_{a_1-1}}{\a_j-\a_{j+1}}=5.5$ for $j=1$ and 
\item $\frac{\l'_1- \l'_{a_j}+ \l'_{a_j+1}+\m'_{a_1}-\m'_{a_1-1}}{\a_j-\a_{j+1}}=8$ for $j=2$.
\end{itemize} The following LR tableaux are of shape $\A(1)$ resp. $\A(2)$ and contain fewer than $\a_k$ entries $k$ in row $a_k$, i.e. fewer than one entry $2$ in row $4$.

\[\D^1=\young(:::::::::::1111,:::::::::12222,::::::1123333,:::1134444,:1122,:22,:33,1)\]
\[\D^2=\young(::::::::::::::1111,::::::::::::12222,::::::1111123333,:::11134444,:1122,:22,:33,1)\]

For $j=1$ the following LR tableaux are the $C_n$ from the above construction for $n=3,4,5$.

\[C_3=\young(:::::::::::::::::1111,:::::::::::::::11222,::::::1111111222333,:::222233444,:2233,:33,:44,2)\]
\[C_4=\young(::::::::::::::::::::1111,::::::::::::::::::11222,::::::1111111111222333,:::2222222334,:2233,:33,:44,2)\]
\[C_5=\young(:::::::::::::::::::::::1111,:::::::::::::::::::::11222,::::::1111111111111222333,:::22222222223,:2233,:33,:44,2)\]

For $j=2$ the following LR tableaux are the $C_n$ from the above construction for $n=5,6$ (there are also $C_n$ for $n=3,4,7,8$ which we do not present here).

\[C_5=\young(:::::::::::::::::::::::1111,:::::::::::::::::::::11222,::::::1111111111111222333,:::11122222333,:3333,:44,:55,3)\]
\[C_6=\young(::::::::::::::::::::::::::1111,::::::::::::::::::::::::11222,::::::1111111111111111222333,:::111222223334,:3333,:44,:55,3)\]

The above constructions prove $Q^{\l,\m}_{\l',\m'}(m)>Q^{\l,\m}_{\l',\m'}(m-1)>\ldots>Q^{\l,\m}_{\l',\m'}(0)$ in the case $\l=(\a_1^{a_1},\a_2,\a_3,\ldots\a_k), \m=(\a_1^{a_1-1})$.

In the more general case there can be $i$ with $\m_i=\l_i>\l_{i+1}$ and $\m'_i<\l'_{i+1}$ (so the rows $i$ and $i+1$ of $\A(0)=\l'/\m'$ are connected). We notice that for $n< \frac{ \l'_{i+1}-\m'_i}{\m_i-\l_{i+1}}$ we can construct an LR tableau $C_n$ of shape $\A(n)$ such that row $i+1$ contains $\m'_i-\m'_{i+1}+n(\m_i-\m_{i+1})$ times the entry $1$. Furthermore, we notice that no LR tableau of shape $\A(n-1)$ can contain $\m'_i-\m'_{i+1}+n(\m_i-\m_{i+1}) - (\l_{i+1}-\m_{i+1})$ entries $1$ in  row $i+1$ because there are not enough boxes in row $i+1$ without a box directly on top. So we again have $Q(n)>Q(n-1)$ for these $n$ and for the other $n$ we can specialise to the above case with $\l=(\a_1^{a_1},\a_2,\a_3,\ldots\a_k),\: \m=(\a_1^{a_1-1})$.
\end{proof}

\begin{Bsp}
 Let $\l'=(7^2,5,4^3,3,2^2),\m'=(4,3^3,2),\l=(1^5),\m=(1^2)$. So
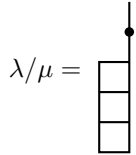
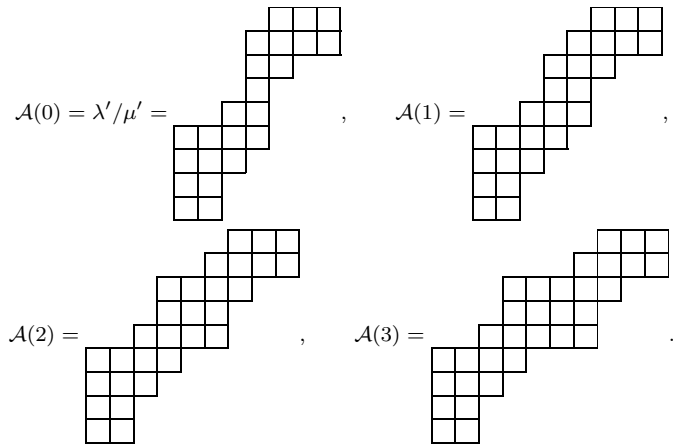
\begin{figure}[h]
\psset{xunit=0.4cm,yunit=0.4cm,runit=0.4cm}
\begin{pspicture*}(-4.2,-0.1)(1.2,5)
\put(-3,2.5){$\lm=$}\psline(1,5)(1,0)(0,0)(0,3)(1,3) \psline(0,2)(1,2)\psline(0,1)(1,1) \psdot(1,4)
\end{pspicture*}
\end{figure}

and
{\footnotesize
\[\A(0)=\l'/\m'=\young(::::\h\h\h,:::\h\h\h\h,:::\h\h,:::\h,::\h\h,\h\h\h\h,\h\h\h,\h\h,\h\h),\qquad \A(1)=\young(:::::\h\h\h,::::\h\h\h\h,:::\h\h\h,:::\h\h,::\h\h\h,\h\h\h\h,\h\h\h,\h\h,\h\h),\]
\[ \A(2)=\young(::::::\h\h\h,:::::\h\h\h\h,:::\h\h\h\h,:::\h\h\h,::\h\h\h\h,\h\h\h\h,\h\h\h,\h\h,\h\h), \qquad
\A(3)=\young(:::::::\h\h\h,::::::\h\h\h\h,:::\h\h\h\h\h,:::\h\h\h\h,::\h\h\h\h\h,\h\h\h\h,\h\h\h,\h\h,\h\h) .\]}
By Theorem~\ref{Sa:main}, we have for $n\geq m=7$ that $Q(n)=Q(7)>Q(6)>\ldots>Q(0)$.
In fact, we have
\begin{align*}
Q(0) && Q(1)  && Q(2)  && Q(3)   && Q(4)   && Q(5)   && Q(6)   && Q(n\geq7)\\
2184 && 26.421 && 92.030 && 172.795 && 229.660 && 254.420 && 260.761 && 261.512.
\end{align*}
\end{Bsp}

\begin{Bsp}
Let $\l=(6,5,3,2,1),\m=(6,1^4),\l'=(8^2,5,3^2,2,1)$ and $\m'=(4,3,2,1^2)$. So
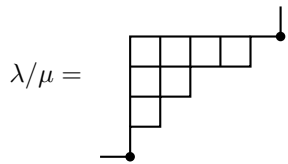
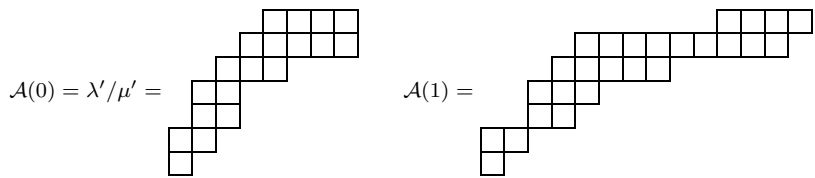
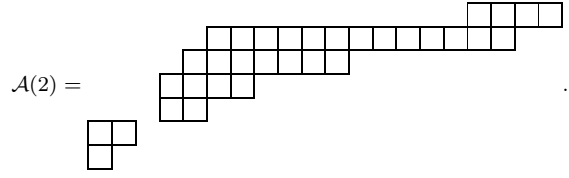
\begin{figure}[h]
\psset{xunit=0.4cm,yunit=0.4cm,runit=0.4cm}
\begin{pspicture*}(-4.2,-0.2)(6.2,5.2)
\put(-3,2.5){$\lm=$}
\psline(6,5)(6,4)(5,4)(5,3)(3,3)(3,2)(2,2)(2,1)(1,1)(1,0)(0,0)\psdot(6,4)\psdot(1,0)
\psline(5,4)(1,4)(1,1) \psline(4,4)(4,3)\psline(3,4)(3,3)\psline(2,4)(2,2)\psline(1,3)(3,3) \psline(1,2)(2,2)
\end{pspicture*}
\end{figure}

and
{\footnotesize
\[\A(0)=\l'/\m'=\young(::::\h\h\h\h,:::\h\h\h\h\h,::\h\h\h,:\h\h,:\h\h,\h\h,\h) \qquad \A(1)=\young(::::::::::\h\h\h\h,::::\h\h\h\h\h\h\h\h\h,:::\h\h\h\h\h,::\h\h\h,::\h\h,\h\h,\h)\]
\[\A(2)=\young(::::::::::::::::\h\h\h\h,:::::\h\h\h\h\h\h\h\h\h\h\h\h\h,::::\h\h\h\h\h\h\h,:::\h\h\h\h,:::\h\h,\h\h,\h).\]
}
By Theorem~\ref{Sa:main}, there exists an $m$ with $Q(n)=Q(m)$ for $n\geq m$ but we cannot use the given formula. For $n=0$ the skew diagram $\A(n)$ is connected, for $1\leq n <4$ $\A(n)$ decays into two skew diagrams, one containing the top  five rows and one the rows below row $5$. For $4 \leq n$ the skew diagram decays into three skew diagrams, one containing the topmost row, one containing the rows $2$ to $5$ and one containing the rows below. Deleting the empty columns in $\A(4)$ and ignoring the empty columns of $\lm$ which only translate the disconnected skew diagrams we can now use the formula on $\widetilde{\A(4)}=(29,25,14,8,4,2,1)/(25,4,3,2,2)$ and $\widetilde{\lm}=(4,4,2,1)/(4)$ which gives $\widetilde{m}=4$.
So in total we have for $n\geq m=8=4+\widetilde{m}$ that $Q(n)=Q(8)>Q(7)>\ldots>Q(0)$.
In fact, we have
\begin{align*}
Q(0) && Q(1)  && Q(2)  && Q(3)   && Q(4)   && Q(5)   && Q(6)   && Q(7)   && Q(n\geq8)\\
910  && 18\,271 && 38\,016 && 49\,635  && 54\,176  && 55\,480  && 55\,826  && 55\,889  && 55\,895.
\end{align*}
\end{Bsp}

\section{The Behaviour of $P^{\l,\m,\n}_{\l',\m',\n'}(n)$} \label{sec:P}

For $c(\l;\m,\n),c(\l';\m',\n')\neq 0$ we define $P^{\l,\m,\n}_{\l',\m',\n'}(n)=c(n\l+\l';n\m+\m',n\n+\n')$ and write simply $P(n)$ if $\l,\m,\n,\l',\m',\n'$ are known from the context.

\begin{Le} \label{Le:P}
 Let  $c(\l';\m',\n')>0$. Let $\lm, \l/\n$ or $\left((\l_1)^{l(\l)}/\m \right)^\circ/\n$ be a partition or a rotated partition.
Then there exists an integer $m$ with
\begin{align*}
 P^{\l,\m,\n}_{\l',\m',\n'}(n)&=P^{\l,\m,\n}_{\l',\m',\n'}(m)&& \textnormal{ for }n\geq m.
\end{align*}
 
\end{Le}
\begin{proof}
 This follows directly from Theorem~\ref{Sa:main}. In the case that $\left((\l_1)^{l(\l)}/\m \right)^\circ/\n$ is a partition we have to use rotation symmetry and $c(\l;\m,\n)=c(\l;\n,\m)$.
\end{proof}

\begin{Bem}
Note that the conditions of Lemma~\ref{Le:P} force  $c(\l;\m,\n)=1$ or  $c(\l;\m,\n)=0$.

 Note furthermore, that we can use the formula in Theorem~\ref{Sa:main} to obtain an $m$ with $P(n)=P(m)$ for $n\geq m$ but the $m$ obtained by the formula in Theorem~\ref{Sa:main} does not have to be minimal.
\end{Bem}

\begin{Le} \label{Le:Rassartg}
 Let $c(\l';\m',\n')>0$.
Then there exist an integer $m$ and a polynomial $g(n)$ with 
\begin{align*}
 P^{\l,\m,\n}_{\l',\m',\n'}(n)&=g(n)&& \textnormal{ for }n\geq m.
\end{align*}
\end{Le}

\begin{proof}
 This follows directly from the work of Rassart \cite{Rassart} mentioned in the introduction. Let $k=\max(l(\l),l(\m),l(\n),l(\l'),l(\m'),l(\n'))$ be the maximal length of the partitions involved. The LR chamber complex $\LRk\subseteq\Rk$ contains those triples of partitions $(\a,\b,\c)$ which have positive LR coefficient $c(\a;\b,\c)$. This chamber complex decays into cones in which the LR coefficient of the triple $(\a,\b,\c)$ is given by a polynomial in the $3k$ variables $\a_1,\ldots,\a_k,\b_1,\ldots,\c_k$. The LR coefficients of triples which lie on a wall between two cones are also given by a polynomial of those variables.

From this it follows that the stretched LR coefficient $c(n\l;n\m,n\n)$ for a fixed triple of partitions $(\l,\m,\n)$ is given by a polynomial in $n$. Suppose $(\l,\m,\n)$ lies inside a cone whose LR coefficients are given by the polynomial $r(\l_1,\ldots,\n_k)$. Since the stretched triple $(n\l,n\m,n\n)$ lies inside the same cone, these LR coefficients are given by $r(n\l_1,\ldots,n\n_k)$, which is a polynomial in $n$ for fixed partitions $\l,\m,\n$. The same applies if $(\l,\m,\n)$ lies not inside a cone but instead on a wall, since then $(n\l,n\m,n\n)$ will lie on the same wall and is therefore given by the same polynomial.

Let us now look at the generalised stretched LR coefficients $P(n)=P^{\l,\m,\n}_{\l',\m',\n'}(n)=c(n\l+\l';n\m+\m',n\n+\n')$. Assume that the triple $(\l,\m,\n)$ lies inside a cone, in which the LR coefficients are given by the polynomial $r(\l_1,\ldots,\n_k)$. Now $(\l',\m',\n')$ may lie in another cone, as may  $(\l+\l',\m+\m',\n+\n')$ and  $(2\l+\l',2\m+\m',2\n+\n')$ and so on. But the lines $\{ (n\l,n\m,n\n)| n\in\mathbb{N}\}$ and $\{ (n\l+\l',n\m+\m',n\n+\n')| n\in\mathbb{N}\}$ are parallel. So for large $n$ the triple $(n\l+\l',n\m+\m',n\n+\n')$ has to lie in the same cone as the triple $(\l,\m,\n)$. Therefore, $P(n)$ is given for large $n$ by the polynomial $r(n\l_1+\l_1',\ldots,n\n_k+\n_k')$ which is a polynomial in $n$ for fixed partitions. 

Now suppose  that  the triple $(\l,\m,\n)$ lies on a wall between two cones. If the triple $(\l',\m',\n')$ lies on the same wall the same argument as above applies. If the triple  $(\l,\m,\n)$ lies in a cone then the triple  $(n\l+\l',n\m+\m',n\n+\n')$ will, for large $n$, also lie in a fixed cone, and by the same argument as above $P(n)$ will be given by a polynomial for large $n$.
\end{proof}

\begin{Bsp}
 Let $\l=(6,5,4,3^2,1),\m=(5,3,2,1),\n=(5,3,2,1)$. We then have $c(\l;\m,\n)=12$ and the polynomial
\[c(n\l;n\m,n\n)=\frac{(n+1)(n+2)(n+3)(n+4)(n+5)(2n^2+5n+7)}{840} \] is of degree $7$.

Let $\l'=(9^3,7,3^4,2,1),\m'=(7^2,3,2^3,1^2),\n'=(8,5,3^2,2^2,1)$. 
We then have $c(\l';\m',\n')=39$ and
\[
\begin{array}{c|c|c|c|c|c|c}
n:    &       0&       1&        2&          3&           4& n\geq5 \\
\hline
P(n): & 39     & 30\,920 & 509\,202 & 3\,101\,626 & 12\,098\,348 & g(n)   \\
g(n): & 55\,407 & 50\,333 & 513\,782 & 3\,102\,223 & 12\,098\,382 & g(n)
\end{array}\]
with
\begin{align*}
 g(n)=\frac{1}{360}\bigl(&8490n^7+214\,525n^6+1\,664\,232n^5+5\,835\,910n^4+904\,140n^3\\&+8\,621\,725n^2-19\,075\,662n+19\,946\,520\bigr).
\end{align*} (We checked $P(n)=g(n)$ for $5\leq n \leq 17$ by computer. It is still possible that $(n\l+\l',n\m+\m',n\n+\n')$ moves to another cone for higher $n$.)
The generating function $G(z)=\sum_n g(n)z^n$ of $g(n)$ is given by:
\begin{align*}
  G(z)=\frac{1}{(1-z)^8}\bigl(& -141\,993z^7+752\,295z^6-1\,841\,275z^5+2\,726\,336z^4\\&-2\,701\,501z^3+1\,662\,514z^2-392\,923z+55\,407\bigr).
\end{align*}

\end{Bsp}

Many calculations suggest that Lemma~\ref{Le:P} and Lemma~\ref{Le:Rassartg} can be generalised.

\begin{Con} \label{Conjecture-P}
Let $f(n)=c(n\l;n\m,n\n)$ be a polynomial of degree $d$. Let $c(\l';\m',\n')\neq 0$. Then there exist a polynomial $g(n)$ of degree $d$ and an integer $m$ such that $P^{\l,\m,\n}_{\l',\m',\n'}(n)=g(n)$ for $n\geq m$.

In particular for $c(\l;\m,\n)=1$ there exists an integer $m$ with $P(n)=P(m)$ for $n\geq m$.
\end{Con}

We will say that a triple of partitions $(\l,\m,\n)$ is larger than another triple $(\l',\m',\n')$ if there exist triples $(\l^i,\m^i,\n^i)$ with $c(\l^i;\m^i,\n^i)\neq0$ with
\begin{align*}
 \l&=\left(\cdots\left(\left(\l'+\l^1\right)+\l^2\right)\cdots\right)+\l^n,\\
\m&=\left(\cdots\left(\left(\m'+\m^1\right)+\m^2\right)\cdots\right)+\m^n,\\
\n&=\left(\cdots\left(\left(\n'+\n^1\right)+\n^2\right)\cdots\right)+\n^n.
\end{align*}

Since the $+$ operation is commutative $(\l,\m,\n)$ is larger than $(\l',\m',\n')$ if and only if $c(\l-\l';\m-\m',\n-\n')> 0$.

\begin{Le}\label{Le:P2}
Let $f(n)=c(n\l;n\m,n\n)$ be a polynomial of degree $d$. Let a multiple of the triple $(\l,\m,\n)$ be larger than the triple $(\l',\m',\n')$.
Then there exists a polynomial $g(n)$ of degree $d$ and an integer $m$ such that $P^{\l,\m,\n}_{\l',\m',\n'}(n)=g(n)$ for $n\geq m$. 
\end{Le}
\begin{proof}
Choose $k$ such that $(k\l,k\m,k\n)$ is larger than $(\l',\m',\n')$.

By Lemma~\ref{Le:Rassartg}, there exist a polynomial $g(n)$ and an integer $m$ such that $P(n)=g(n)$ for $n\geq m$. Suppose in the following that $n\geq m$. We now have $g(n)\geq f(n)$ by Lemma~\ref{Le:add}. But since $(k\l,k\m,k\n)$ is larger than $(\l',\m',\n')$ we also have $f(k+n)\geq g(n)$, also by Lemma~\ref{Le:add}. Since both $f(n)$ and $f(k+n)$ have degree $d$ and $f(k+n)\geq g(n)\geq f(n)$ it follows that $g(n)$ has to be of degree $d$ also.
\end{proof}

{\bfseries Acknowledgement:} John Stembridge's ``SF-package for maple'' \cite{stemmaple} and A. S. Buch's ``Littlewood-Richardson Calculator" \cite{LRcalc} were very helpful for computing examples. Furthermore, my thanks go to Etienne Rassart, Emmanuel Briand, Christine Bessenrodt, Ron King and Martin Rubey for helpful discussions. This paper was inspired by a talk given by Ron King at SLC 60 about stretched LR coefficients.

\end{document}